\documentclass[11pt]{article}
\usepackage{graphicx,amsmath,amsbsy,amssymb,citesort,epsfig}
\usepackage{algorithm,algorithmic}

\def\real    { \mathbb{R} }
\def\reals   { \mathbb{R} }

\newcommand{\qed}{{\unskip\nobreak\hfil\penalty50\hskip2em\vadjust{}
           \nobreak\hfil$\Box$\parfillskip=0pt\finalhyphendemerits=0\par}}
\newtheorem{thm}{Theorem}[section] 
\newtheorem{lemma}{Lemma}[section] 
\newtheorem{cor}{Corollary}[section]

\newcommand{\bitem}{\begin{itemize}}
\newcommand{\eitem}{\end{itemize}}

\newcommand{\argmax}{\mathrm{argmax}}
\newcommand{\argmin}{\mathrm{argmin}}

\newcommand{\supp}{\mathrm{supp}}
\newcommand{\beqn}{\begin{equation}}
\newcommand{\eeqn}{\end{equation}}
\newcommand{\balign}{\begin{align}}
\newcommand{\ealign}{\end{align}}

\newcommand{\norm}[1]{\| #1 \|}

\def \Ran {{\mathcal{R}}}

\title{{ Analysis of Orthogonal Matching Pursuit using the \\ Restricted Isometry Property}}

\author{Mark~A.~Davenport and~Michael~B.~Wakin
\thanks{M. Davenport is with the Department of Electrical and Computer Engineering,
Rice University, Houston, TX 77098, USA. e-mail: md@rice.edu. M.
Wakin is with Division of Engineering, Colorado School of Mines,
Golden, CO 80401, USA. e-mail: mwakin@mines.edu. This research was
partially supported by NSF Grant CCF-0830320, DARPA Grants
N66001-08-1-2065 and HR0011-08-1-0078, AFOSR Grants FA9550-07-1-0301
and FA9550-09-1-0465, and ONR Grant N00014-07-1-0936. Thanks to
Marco Duarte and Chinmay Hegde for their comments on a preliminary
version of this manuscript.}
}

\date{August 2009}

\begin{document}

\maketitle

\vspace{-0.3in}

\begin{abstract}
Orthogonal Matching Pursuit (OMP) is the canonical greedy algorithm
for sparse approximation.  In this paper we demonstrate that the
restricted isometry property (RIP) can be used for a very
straightforward analysis of OMP. Our main conclusion is that the RIP
of order $K+1$ (with isometry constant $\delta <
\frac{1}{3\sqrt{K}}$)
is sufficient for OMP to exactly recover any $K$-sparse signal. Our
analysis relies on simple and intuitive observations about OMP and
matrices which satisfy the RIP. For restricted classes of $K$-sparse
signals (those that are highly compressible), a relaxed bound on the
isometry constant is also established. A deeper understanding of
OMP may benefit the analysis of greedy algorithms in
general. To demonstrate this, we also briefly revisit the analysis
of the Regularized OMP (ROMP) algorithm.
\end{abstract}

\section{Introduction}

\subsection{Orthogonal Matching Pursuit}

Orthogonal Matching Pursuit (OMP) is the canonical greedy algorithm
for sparse approximation. Letting $\Phi$ denote a matrix of size $M
\times N$ (where typically $M < N$) and $y$ denote a vector in
$\real^M$, the goal of
OMP is to recover a coefficient
vector $\widehat{x} \in \real^N$ with roughly $K < M$ nonzero terms
so that $\Phi \widehat{x}$ equals $y$ exactly or approximately. OMP
is frequently used to find sparse representations for signals $y \in
\real^M$ in settings where $\Phi$ represents an overcomplete
dictionary for the signal space~\cite{PatiOMP,DavisOMP,Tropp04}. It
is also commonly used in compressive sensing (CS), where $y = \Phi
x$ represents compressive measurements of a sparse or nearly-sparse
signal $x \in \real^N$ to be
recovered~\cite{DonohoCS,CandesCS,TroppOMP}.

One of the attractive features of OMP is its simplicity. The entire
algorithm is specified in Algorithm~\ref{alg:omp}, and it requires
approximately the same number of lines of code to implement in a
software package such as Matlab. Despite its simplicity,
OMP is empirically competitive in terms of approximation
performance~\cite{Tropp04,rauhut2}.

Theoretical analysis of OMP to date has concentrated primarily on
two fronts. The first has involved the notion of a coherence
parameter $\mu := \max_{i,j} |\langle \phi_i, \phi_j \rangle|$,
where $\phi_i$ denotes column $i$ of the matrix $\Phi$. When the
columns of $\Phi$ have unit norm and $\mu < \frac{1}{2K-1}$, it has
been shown~\cite{Tropp04} that OMP will recover any $K$-sparse
signal $x$ from the measurements $y = \Phi x$. This guarantee is
deterministic and applies to any matrix $\Phi$ having normalized
columns and $\mu < \frac{1}{2K-1}$.

The second analytical front has involved the notion of probability.
Suppose $x \in \real^N$ with $\|x\|_0 := |\supp(x)| \le K$ and that $\Phi$ is
drawn from a suitable random distribution (independently of $x$)
with $M = O (K \log(N))$ rows. Then with high probability,
OMP will recover $x$ exactly from the measurements $y = \Phi
x$~\cite{TroppOMP}. It is not guaranteed, however, that any such
fixed matrix will allow recovery of all sparse $x$ simultaneously.

%

\begin{algorithm}[t]
\caption{Orthogonal Matching Pursuit} \label{alg:omp}
\begin{algorithmic}
\STATE \textbf{input:} $\Phi$, $y$, stopping criterion
\STATE \textbf{initialize:} $r^0 = y$, $x^0 = 0$, $\Lambda^0 =
\emptyset$, $\ell = 0$
\WHILE{not converged}
\STATE
\begin{tabular}{ll}
\textbf{match:} & $h^{\ell} = \Phi^T r^{\ell}$ \\
\textbf{identify:} & $\Lambda^{\ell+1} = \Lambda^\ell \cup
\{\arg\max_j |h^\ell(j)| \}$ \; (if multiple maxima exist, choose
only one)
\\
\textbf{update:} & $x^{\ell+1} = \argmin_{z: ~ \supp(z) \subseteq
\Lambda^{\ell+1}}  \|y - \Phi z\|_2$ \\
 & $r^{\ell+1} = y - \Phi x^{\ell+1}$   \\
 & $\ell = \ell+1$
\end{tabular}
\ENDWHILE
\STATE \textbf{output:} $\widehat{x} = x^{\ell} =\argmin_{z: ~
\supp(z) \subseteq \Lambda^{\ell}} \|y - \Phi z\|_2$
\end{algorithmic}
\end{algorithm}

\subsection{The Restricted Isometry Property}

As an alternative to coherence and to probabilistic analysis,
a large number of algorithms within the broader field of CS have
been studied
using the {\em restricted isometry property}
(RIP) for the matrix $\Phi$~\cite{CandesDLP}. A matrix $\Phi$
satisfies the RIP of order $K$ if there exists a constant $\delta
\in (0,1)$ such that
\begin{equation} \label{eq:RIP}
(1-\delta)\|x\|_2^2 \le \|\Phi x \|^2_2 \le (1+\delta)\|x\|_2^2
\end{equation}
holds for all $x$ such that $\|x\|_0 \le K$. In other words, $\Phi$
acts as an approximate isometry on the set of vectors that are
$K$-sparse. Much is known about finding matrices that satisfy the
RIP. For example, if we draw a random $M \times N$ matrix $\Phi$
whose entries $\phi_{ij}$ are independent and identically
distributed sub-Gaussian random variables, then provided that
\begin{equation}
M = O\left(\frac{K\log(N/K)}{\delta^2}\right), \label{eq:ripm}
\end{equation}
with high probability $\Phi$ will satisfy the RIP of order
$K$~\cite{jlcs,DeVoreL1IO}.

When it is satisfied, the RIP for a matrix $\Phi$
provides a
sufficient condition to guarantee successful sparse recovery using a
wide variety of
algorithms~\cite{CandesDLP,CandesRIP,MIP_cai,romp,romp2,cosamp,sp,Thresh,IHT,chartrand2008restricted}.
As an example, the RIP of order $2K$ (with isometry constant
$\delta < \sqrt{2}-1$)
is
a sufficient condition
to permit $\ell_1$-minimization (the canonical convex optimization
problem for sparse approximation) to exactly recover any $K$-sparse
signal and to approximately recover those that are nearly
sparse~\cite{CandesRIP}. The same RIP assumption is also a
sufficient condition for robust recovery in noise using a modified
$\ell_1$-minimization~\cite{CandesRIP}.

Despite the considerable attention that has been paid to both OMP
and the RIP, analysis of OMP using the RIP has been relatively
elusive to date. However, several alternative greedy algorithms have
been proposed---all essentially modifications of OMP---that are
apparently much more amenable to RIP-based analysis. The Regularized
Orthogonal Matching Pursuit (ROMP)~\cite{romp,romp2} and Subspace
Pursuit (SP)~\cite{sp} algorithms differ from OMP in the
identification step, while the Compressive Sampling Matching Pursuit
(CoSaMP)~\cite{cosamp} and DThresh~\cite{Thresh} algorithms differ
from OMP in both the identification and the
update steps.  For each of these algorithms it has been shown that
the RIP of order $C K$ (where $C \ge 2$ is a constant depending on
the algorithm) with $\delta$ adequately small is sufficient for
exact recovery of $K$ sparse signals.

\subsection{Contributions}

Our contributions in this paper are twofold. First, we begin in
Section~\ref{sec:obs} with some very simple observations regarding
OMP. Many of these facts are known to practitioners in
the field but may not be obvious to a novice, and we feel that such
readers may find value in a short exposition.

Critically, these observations also set the stage for our main
results in Section~\ref{sec:analysis}, in which we demonstrate that
the RIP can be used for a very straightforward analysis of OMP. Our
analysis revolves around three key facts: (1)~that in each step of
the algorithm, the residual vector $r^\ell$ can be written as a
matrix times a sparse signal, (2)~that this matrix satisfies the
RIP, and (3)~that consequently a sharp bound can be established for
the vector $h^\ell$ of inner products. Our main conclusion,
Theorem~\ref{theo:main}, states that the RIP of order
$K+1$ (with $\delta < \frac{1}{3\sqrt{K}}$) is sufficient for OMP to
exactly recover any $K$-sparse signal in exactly $K$ iterations.
However, for restricted classes of $K$-sparse signals (those with
sufficiently strong decay in the nonzero coefficients),
a relaxed bound on the isometry constant can be used. We discuss
such extensions of our results in Section~\ref{sec:extensions}. A
deeper understanding of
OMP may also
benefit the analysis
of greedy algorithms in general. To demonstrate this, we
briefly revisit the analysis of the ROMP algorithm in
Section~\ref{sec:extensions}.

%

\subsection{Context}

Let us place
Theorem~\ref{theo:main}
in the context of the OMP
literature. Using the RIP as a sufficient condition to guarantee OMP
performance is apparently novel. Moreover, the fact that our bound
requires only the RIP of order $K+1$ is apparently unique among the
published CS literature; much more common are results requiring the
RIP of order $1.75K$~\cite{MIP_cai}, $2K$~\cite{CandesRIP,romp},
$3K$~\cite{sp,IHT}, $4K$~\cite{cosamp}, and so on. Of course, such
results
often permit the isometry constant to be much
larger.\footnote{Note that a smaller order of the
RIP is not necessarily a weaker requirement if the required constant
is also significantly smaller. For example, Corollary 3.4
of~\cite{cosamp} implies that if $\Phi$ satisfies the RIP of order
$K+1$ with constant $\delta$, then $\Phi$ also satisfies the RIP of
order $2K$ with constant $4 \delta $.}

If one wishes to use the RIP of order $K+1$ as a sufficient
condition for exact recovery of all $K$-sparse signals via OMP (as
we have),
then little improvement is possible in relaxing the isometry
constant $\delta$ above $\frac{1}{3\sqrt{K}}$. In particular, there
exists a matrix satisfying the RIP of order $K+1$ with $\delta \le
\frac{1}{\sqrt{K}}$ for which there exists a $K$-sparse signal $x
\in \real^N$ that cannot be recovered exactly via $K$ iterations of
OMP. (This is conjectured in~\cite{sp} with a suggestion for
constructing such a matrix, and for the case $K=2$ we have confirmed
this via experimentation.)


Unfortunately, from (\ref{eq:ripm}) we see that finding a matrix
$\Phi$ satisfying the RIP of order $K+1$ with an isometry constant
$\delta < \frac{1}{3\sqrt{K}}$ may require $M = O(K^2 \log(N/K))$
random measurements. If one wishes to guarantee exact recovery of
all $K$-sparse signals via OMP (as we have),
then
little improvement is
possible in relaxing this number. In particular, it has been
argued~\cite{rauhut2008impossibility} that when $M \lesssim
K^{3/2}$, for most random $M \times N$ matrices $\Phi$ there will
exist some $K$-sparse signal $x \in \real^N$ that cannot be
recovered exactly via $K$ iterations of OMP.

It is also worth comparing our RIP-based analysis with
coherence-based analysis~\cite{Tropp04}, as both techniques provide
a sufficient condition for OMP to recover all $K$-sparse signals. It
has been shown~\cite{TroppOMP} that in a random $M \times N$ matrix,
the coherence parameter $\mu$ is unlikely to be smaller than
$\log(N)/\sqrt{M}$. Thus, to ensure $\mu < \frac{1}{2K-1}$, one
requires $M = O(K^2 \log^2(N))$, which is roughly the same as what
is required by our analysis. We note that neither result is strictly
stronger than the other; we have confirmed experimentally that there
exist matrices that satisfy our RIP condition but not the coherence
condition, and vice versa.

Finally, we note that the aforementioned modifications of OMP (the
ROMP, SP, CoSaMP, and DThresh algorithms) all have RIP-based
guarantees of robust recovery in noise and stable recovery of
non-sparse signals. To date, no such RIP-based or coherence-based
guarantees have been provided for OMP itself. We speculate that our
perspective may help to further the understanding of OMP and perhaps
provide a route to such a guarantee. At present, however, this
remains a topic of ongoing work.

%
%
%
%
%
%

\subsection{Notation}

Before proceeding, we set our notation.
Suppose $\Lambda \subset \{1,2,\dots,N\}$. We let $\Lambda^c =
\{1,2,\dots,N\} \backslash \Lambda$. By $x|_\Lambda$ we mean the
length $|\Lambda|$ vector containing the entries of $x$ indexed by
$\Lambda$.

By $\Phi_\Lambda$ we mean the $M \times |\Lambda|$ matrix obtained
by selecting the columns of $\Phi$ indexed by $\Lambda$, and by
$\Ran(\Phi_\Lambda)$ we mean the range, or column space, of
$\Phi_\Lambda$. We will assume throughout that when $|\Lambda| \le
M$, $\Phi_\Lambda$ is full rank, in which case we let
$\Phi_\Lambda^\dag := (\Phi_{\Lambda}^T \Phi_{\Lambda})^{-1}
\Phi_{\Lambda}^T$ denote the Moore-Penrose pseudoinverse of
$\Phi_\Lambda$.


We denote the orthogonal projection operator onto
$\Ran(\Phi_{\Lambda})$ by $P_{\Lambda} := \Phi_{\Lambda}
\Phi_{\Lambda}^\dag$. Similarly, $P^\perp_{\Lambda} := (I -
P_{\Lambda})$ is the orthogonal projection operator onto the
orthogonal complement of $\Ran(\Phi_{\Lambda})$. We note that any
orthogonal projection operator $P$ obeys $P = P^T = P^2$.

Finally, we define $A_{\Lambda} := P^\perp_{\Lambda} \Phi$. This
matrix is the result of orthogonalizing the columns of $\Phi$
against $\Ran(\Phi_{\Lambda})$. It is therefore equal to zero on
columns indexed by $\Lambda$.


\section{Observations}
\label{sec:obs}



Let us begin with some very simple observations regarding
OMP as presented in Algorithm~\ref{alg:omp}. The key idea is
to try to iteratively estimate a set $\Lambda$ that contains the
locations of the
nonzeros
of $x$ by starting with $\Lambda =
\emptyset$ and then adding a new element to $\Lambda$ in each
iteration. In order to select which element to add, the algorithm
also maintains a residual vector $r \notin \Ran(\Phi_\Lambda)$ that
represents
the component of the measurement vector $y$ that cannot be explained
by the columns of $\Phi_\Lambda$.
Specifically, at the beginning of the $\ell^{\mathrm{th}}$
iteration, $\Lambda^\ell$ is our current estimate of $\supp(x)$, and
the residual $r^\ell$ is defined as $r^\ell = y - \Phi x^\ell$ where
$\supp(x^\ell) \subseteq \Lambda^\ell$. The element added to
$\Lambda^\ell$ is the index of the column of $\Phi$ that has the
largest inner product with $r^\ell$.

Our first observation is that $r^\ell$ can be viewed as the
orthogonalization of $y$ against the previously chosen columns of
$\Phi$. To see this, note that the solution to the least squares
problem in the update step is given by
\begin{equation}
x^{\ell}|_{\Lambda^{\ell}} = \Phi_{\Lambda^{\ell}}^\dag y
\mathrm{~~~and~~~} x^{\ell}|_{(\Lambda^{\ell})^c} = 0. \label{eq:ls}
\end{equation}
Thus we observe that
$$
r^{\ell} = y - \Phi x^{\ell}
 = y - \Phi_{\Lambda^{\ell}} \Phi_{\Lambda^{\ell}}^\dag y
 = (I - P_{\Lambda^{\ell}}) y
 =  P^\perp_{\Lambda^{\ell}} y.
$$
Note that it is not actually necessary to explicitly compute
$x^{\ell}$ in order to calculate $r^\ell$.

Our second observation is that, in the matching step, one may
correlate $r^\ell$ either with the columns of $\Phi$ or with the
columns of $A_{\Lambda^\ell}$. To see this equivalence, observe that
$ r^\ell = P^\perp_{\Lambda^{\ell}} y = P^\perp_{\Lambda^{\ell}}
P^\perp_{\Lambda^{\ell}} y = (P^\perp_{\Lambda^{\ell}})^T
P^\perp_{\Lambda^{\ell}} y $
and so
\begin{equation}
h^\ell = \Phi^T r^\ell = \Phi^T (P^\perp_{\Lambda^{\ell}})^T
P^\perp_{\Lambda^{\ell}} y = A_{\Lambda^\ell}^T r^\ell.
\label{eq:har}
\end{equation}
Incidentally, along these same lines
we
observe that
$$
h^\ell  = \Phi^T r^\ell =\Phi^T P^\perp_{\Lambda^{\ell}} y  = \Phi^T
(P^\perp_{\Lambda^{\ell}})^T y = A_{\Lambda^\ell}^T y.
$$
From this we note that it is not actually necessary to explicitly
compute $r^{\ell}$ in order to calculate the inner products during
the matching step; in fact, the original formulation of OMP was stated with instructions
to orthogonalize the remaining columns of $\Phi$ against those
previously chosen
and merely correlate the resulting vectors against
$y$~\cite{PatiOMP,DavisOMP}.
%
%
Additionally, we recall that, in $A_{\Lambda^\ell}$, all columns
indexed by $\Lambda^\ell$ will be zero. It follows that
\begin{equation}
h^\ell(j) = 0 ~~~ \forall j \in \Lambda^\ell,\label{eq:zero}
\end{equation}
and so, since $\Lambda^\ell = \Lambda^{\ell-1} \cup \{j^*\}$ with
$j^* \notin \Lambda^{\ell-1}$,
\begin{equation}
|\Lambda^\ell| = \ell. \label{eq:lamsize}
\end{equation}

Our third observation is that, in the case of noise-free
measurements $y = \Phi x$, we may write
$$
r^\ell =  P^\perp_{\Lambda^{\ell}} y
 = P^\perp_{\Lambda^{\ell}} \Phi x
 = A_{\Lambda^\ell} x.
$$
Again recalling that all columns of $A_{\Lambda^\ell}$ indexed by
$\Lambda^\ell$ are zero, we thus note that when $\supp(x) \subseteq
\Lambda^{\ell}$,  $r^\ell = 0$, and from (\ref{eq:ls}) we also know
that $x^{\ell} = x$ exactly.
It will also be useful to note that for the same reason, we can also
write
\begin{equation}
r^\ell = A_{\Lambda^\ell} \widetilde{x}^\ell, \label{eq:restilde}
\end{equation}
where
\begin{equation}
\widetilde{x}^{\ell}|_{\Lambda^{\ell}} = 0
\mathrm{~~~and~~~} \widetilde{x}^{\ell}|_{(\Lambda^{\ell})^c} = x|_{(\Lambda^{\ell})^c}. \label{eq:xtilde}
\end{equation}


\section{Analysis}
\label{sec:analysis}

Our analysis of OMP will center on the vector $h^\ell$. 
In light of (\ref{eq:har}) and (\ref{eq:restilde}), we
see that $A_{\Lambda^\ell}$ plays a role both in {\em constructing}
and in {\em analyzing} the residual vector.
In Lemma~\ref{lem:PRIP} below, we show that the matrix
$A_{\Lambda^\ell}$ satisfies a modified version of the RIP. This
allows us to very precisely bound the values of the inner products
in the vector $h^\ell$.

We begin with two elementary lemmas whose proofs are given in the
Appendix. Our first result, which is a straightforward
generalization of Lemma 2.1 of~\cite{CandesRIP}, states that RIP
operators approximately preserve inner products between sparse
vectors.

\begin{lemma} \label{lem:ip}
Let $u,v \in \reals^N$ be given, and suppose that a matrix $\Psi$
satisfies the RIP of order $\max(\|u+v\|_0,\|u-v\|_0)$ with isometry
constant $\delta$. Then
\begin{equation} \label{eq:ip-bound}
\left| \langle \Psi u, \Psi v \rangle  - \langle u , v \rangle
\right| \le \delta \|u\|_2 \|v\|_2.
\end{equation}
\end{lemma}

One consequence of this result is that sparse vectors that are
orthogonal in $\reals^N$ remain nearly orthogonal after the
application of $\Psi$. From this observation, it was demonstrated
independently in~\cite{SPARS} and~\cite{sp} that if $\Phi$ has the
RIP, then $A_\Lambda$ satisfies a modified version of the RIP.

\begin{lemma} \label{lem:PRIP}
Suppose that $\Phi$ satisfies the RIP of order $K$ with isometry
constant $\delta$, and let $\Lambda \subset \{1, 2, \ldots, N\}$. If
$|\Lambda| < K$ then
\begin{equation} \label{eq:PRIP}
\left(1 - \frac{\delta}{1-\delta} \right) \|u\|_2^2 \le \| A_\Lambda u\|_2^2 \le (1+\delta) \|u\|_2^2
\end{equation}
for all $u \in \reals^N$ such that $\|u\|_0 \le K - |\Lambda|$ and
$\supp(u) \cap \Lambda =  \emptyset$.
\end{lemma}

In other words, if $\Phi$ satisfies the RIP of order $K$, then
$A_\Lambda$ acts as an approximate isometry on every ($K -
|\Lambda|$)-sparse vector supported on $\Lambda^c$. From
(\ref{eq:restilde}), we recall that the residual vector in OMP is
formed by applying $A_{\Lambda^\ell}$ to a sparse vector supported
on $(\Lambda^\ell)^c$. Combining the above results, then, we may
bound the inner products $h^\ell(j)$ as follows.

\begin{lemma} \label{lem:ip2}
Let $\Lambda \subset \{1, 2, \ldots, N\}$ and suppose $\widetilde{x}
\in \real^N$ with $\supp(\widetilde{x}) \cap \Lambda = \emptyset$.
Define
\begin{equation}
h = A_\Lambda^T A_\Lambda \widetilde{x}. \label{eq:hdef}
\end{equation}
Then if $\Phi$ satisfies the RIP of order $\|\widetilde{x}\|_0 +
|\Lambda| + 1$ with isometry constant $\delta$, we have
\begin{equation}
|h(j) - \widetilde{x}(j)| \le \frac{\delta}{1-\delta}
\norm{\widetilde{x}}_2 \label{eq:hbound}
\end{equation}
for all $j \notin \Lambda$.
\end{lemma}
~~~~~{\em Proof:} From Lemma~\ref{lem:PRIP} we have that the
restriction of $A_\Lambda$ to the columns indexed by $\Lambda^c$
satisfies the RIP of order $ (\|\widetilde{x}\|_0 + |\Lambda| + 1) -
|\Lambda| = \|\widetilde{x}\|_0 + 1$ with isometry constant
$\delta/(1-\delta)$.
%
By the definition of $h$, we also know that
\begin{equation*}
h(j) = \langle A_{\Lambda} \widetilde{x}, A_{\Lambda}
e_j \rangle,
\end{equation*}
where $e_j$ denotes the $j^{\mathrm{th}}$ vector from the cardinal
basis. Now, suppose $j \notin \Lambda$. Then because $\|
\widetilde{x} \pm e_j \|_0 \le \| \widetilde{x}\|_0 + 1$ and
$\supp(\widetilde{x} \pm e_j) \cap \Lambda = \emptyset$, we conclude
from Lemma~\ref{lem:ip} that
$$
| h(j) - \widetilde{x}(j) | = |\langle A_{\Lambda} \widetilde{x},
A_{\Lambda}
e_j \rangle - \langle \widetilde{x}, e_j \rangle|
\le \frac{\delta}{1-\delta} \norm{\widetilde{x}}_2 \norm{e_j}_2.
$$
Noting that $\norm{e_j}_2=1$, we reach the desired conclusion. \qed
\vspace{3mm}

With this bound on the inner products $h^\ell(j)$, we may derive a
sufficient condition under which the identification step of OMP will
succeed.

\begin{cor}\label{cor:j}
Suppose that $\Lambda$, $\Phi$, $\widetilde{x}$ meet the assumptions
specified in Lemma~\ref{lem:ip2}, and let $h$ be as defined in
(\ref{eq:hdef}). If
\begin{equation} \| \widetilde{x}
\|_\infty
> \frac{2\delta}{1-\delta}\norm{\widetilde{x}}_2, \label{eq:infty}
\end{equation}
we are guaranteed that
$
\arg\max_j |h(j)| \in \supp(\widetilde{x}).
$
\end{cor}
~~~~~{\em Proof:} If (\ref{eq:hbound}) is satisfied, then for
indices $j \notin \supp(\widetilde{x})$, we will have $|h(j)| \le
\frac{\delta}{1-\delta} \norm{\widetilde{x}}_2$. (Recall from
(\ref{eq:zero}) that $h(j) =0$ for $j \in \Lambda$.) If
(\ref{eq:infty}) is satisfied, then there exists some $j \in
\supp(\widetilde{x})$ with $ |\widetilde{x}(j)| >
\frac{2\delta}{1-\delta}\norm{\widetilde{x}}_2$. From
(\ref{eq:hbound}) and the triangle inequality, we conclude that for
this index $j$, $|h(j)| >
\frac{\delta}{1-\delta}\norm{\widetilde{x}}_2$. \qed \vspace{3mm}

By choosing $\delta$ small enough, it is possible to guarantee that
the condition (\ref{eq:infty}) is satisfied. In particular, the
lemma below follows from standard arguments.

\begin{lemma} \label{lem:infty2}
For any $u \in \real^N$, $ \|u\|_\infty \ge \|u\|_2/\sqrt{\|u\|_0}.$
%
\end{lemma}

Putting these results together, we can now establish our main
theorem concerning OMP.

\begin{thm}
\label{theo:main} Suppose that $\Phi$ satisfies the RIP of order
$K+1$ with isometry constant $\delta < \frac{1}{3\sqrt{K}}$. Then
for any $x \in \real^N$ with $\|x\|_0 \le K$, OMP will recover $x$
exactly from $y = \Phi x$ in $K$ iterations.
\end{thm}
~~~~~{\em Proof:} The proof works by induction. We start with the
first iteration where $h^0 = \Phi^T \Phi x$ and note that $\Phi =
A_\emptyset$. Because $\|x\|_0 \le K$, Lemma~\ref{lem:infty2} states
that $ \|x\|_\infty \ge \frac{\norm{x}_2}{\sqrt{K}}$. One can also
check
that $\delta < \frac{1}{3\sqrt{K}}$ implies that
$\frac{2\delta}{1-\delta} < \frac{1}{\sqrt{K}}$. Therefore, we are
guaranteed that (\ref{eq:infty}) is satisfied, and so from
Corollary~\ref{cor:j} we conclude that $\arg\max_j |h^0(j)| \in
\supp(x)$.

We now consider the general induction step.  Suppose that we are at
iteration $\ell$ and that all previous iterations have succeeded, by
which we mean that $\Lambda^\ell \subseteq \supp(x)$. From
(\ref{eq:xtilde}), we know that $\supp(\widetilde{x}^\ell) \cap
\Lambda^\ell = \emptyset$ and that $\| \widetilde{x}^\ell\|_0 \le K
- \ell$. From (\ref{eq:lamsize}), we know that $| \Lambda^\ell| =
\ell$. By assumption, $\Phi$ satisfies the RIP of order $K+1 =
(K-\ell) + \ell + 1 \ge \| \widetilde{x}^\ell\|_0 + | \Lambda^\ell|
+ 1$. Finally, using Lemma~\ref{lem:infty2}, we have that
$$
\|\widetilde{x}^\ell\|_\infty \ge
\frac{\norm{\widetilde{x}^\ell}_2}{\sqrt{K-\ell}} \ge
\frac{\norm{\widetilde{x}^\ell}_2}{\sqrt{K}} >
\frac{2\delta}{1-\delta}\norm{\widetilde{x}^\ell}_2.
$$
From Corollary~\ref{cor:j} we conclude that $\arg\max_j
|h^\ell(j)| \in \supp(\widetilde{x}^\ell)$ and hence
$\Lambda^{\ell+1} \subseteq \supp(x)$. \qed \vspace{3mm}

\section{Extensions}
\label{sec:extensions}

\subsection{Strongly-decaying sparse signals}

For even moderate values of the isometry constant $\delta$ there
exist sparse
signals that we can ensure are
recovered exactly. For example, if the decay of coefficients is
sufficiently strong in a sparse signal, we may use
Lemma~\ref{lem:ip2} to ensure that the signal entries are recovered
in the order of their magnitude.

For any $x \in \real^N$ with $\|x\|_0 \le K$ we denote by $x'(j)$
the entries of $x$ ordered by magnitude, i.e.,
$$
|x'(1)| \ge |x'(2)| \ge \cdots \ge |x'(K)| \ge 0
$$
with $x'(K+1) = x'(K+2) = \cdots = x'(N) = 0$.

\begin{thm}
\label{theo:main2} Suppose that $\Phi$ satisfies the RIP of order
$K+1$ with isometry constant $\delta < \frac{1}{3}$. Suppose $x \in
\real^N$ with $\|x\|_0 \le K$ and that for all $j \in
\{1,2,\dots,K-1\}$,
$$
\frac{|x'(j)|}{|x'(j+1)|} \ge \alpha.
$$
If
\begin{equation}
\label{eq:alpha}
 \alpha >
\frac{1+2\frac{\delta}{1-\delta}\sqrt{K-1}}{1-2\frac{\delta}{1-\delta}},
\end{equation}
then OMP will recover $x$ exactly from $y = \Phi x$ in $K$
iterations.
\end{thm}
~~~~~{\em Proof:} The proof again proceeds by induction. At each
stage, OMP will choose the largest entry of $\widetilde{x}^\ell$. To
see this, note that by (\ref{eq:hbound}) we have
$|h^\ell(j) -
\widetilde{x}^\ell(j)| \le \frac{\delta}{1-\delta}
\norm{\widetilde{x}^\ell}_2$. The nonzero entries of
$\widetilde{x}^\ell$ will be comprised of $x'(\ell+1), x'(\ell+2),
\dots, x'(K)$. Thus,
\begin{equation*}
\small \norm{\widetilde{x}^\ell}_2 \le \sqrt{|x'(\ell+1)|^2 +
(K-1)\frac{|x'(\ell+1)|^2}{\alpha^2}} = \frac{|x'(\ell+1)|}{\alpha}
\sqrt{\alpha^2+(K-1)}
\le \frac{|x'(\ell+1)|}{\alpha} (\alpha+\sqrt{K-1}).
\end{equation*}

Now, for the largest entry of $\widetilde{x}^\ell$, we have
\begin{equation}
|h^\ell(j)| \ge |x'(\ell+1)| - \frac{\delta}{1-\delta}
\frac{|x'(\ell+1)|}{\alpha} (\alpha+\sqrt{K-1}) =
\frac{|x'(\ell+1)|}{\alpha} (\alpha - \frac{\delta}{1-\delta}(\alpha
+ \sqrt{K-1})) \label{eq:bound1}
\end{equation}
while for all other entries we have
\begin{equation}
|h^\ell(j)| \le |x'(\ell+2)| + \frac{\delta}{1-\delta}
\frac{|x'(\ell+1)|}{\alpha} (\alpha+\sqrt{K-1}) \le
\frac{|x'(\ell+1)|}{\alpha} (1 +
\frac{\delta}{1-\delta}(\alpha+\sqrt{K-1})) . \label{eq:bound2}
\end{equation}
From (\ref{eq:alpha}), it follows that (\ref{eq:bound1}) is greater
than (\ref{eq:bound2}).
\qed \vspace{3mm}

\subsection{Analysis of other orthogonal greedy algorithms}

We now demonstrate that the techniques used above can also be used
to analyze other orthogonal
greedy algorithms.  We
focus on ROMP~\cite{romp,romp2} for the purpose of
illustration, but similar methods should be able to simplify the
analysis of other orthogonal greedy algorithms such as
SP~\cite{sp}.\footnote{Some of the greedy algorithms that have
been proposed recently, such as CoSaMP~\cite{cosamp} and
DThresh~\cite{Thresh}, do not orthogonalize the residual against the
previously chosen columns at each iteration, and so the techniques
above
cannot be directly applied to these algorithms. However, this
orthogonalization step could easily be added (which in the case of
CoSaMP yields an algorithm nearly identical to
SP).
Orthogonalized versions of these algorithms could then be studied
using these techniques.}

We first briefly describe the difference between ROMP and OMP,
which lies only
in the identification step: whereas OMP adds only
one index to $\Lambda^\ell$ at each iteration, ROMP adds
up to $K$ indices to $\Lambda^\ell$ at each iteration.  Specifically, ROMP
first selects the indices corresponding to the $K$ largest elements
in magnitude of $h^\ell$ (or all nonzero elements of $h^\ell$ if
$h^\ell$ has fewer than $K$ nonzeros), and denotes this set as
$\Omega^{\ell}$.  The next step is to {\em regularize} this set so
that the values are comparable in magnitude.  To do this, define
$R(\Omega^{\ell}) := \{ \Omega \subseteq \Omega^{\ell} : |h^\ell(i)|
\le 2|h^\ell(j)| \ \forall i,j \in \Omega \}$, and set
$$
\Omega^{\ell}_0 := \argmax_{\Omega \in R(\Omega^{\ell})}
\|h^\ell|_\Omega \|_2, 
$$
i.e., $\Omega^{\ell}_0$ is the set with maximal energy among all
regularized subsets of $\Omega^{\ell}$.  Finally, setting
$\Lambda^{\ell+1} = \Lambda^{\ell} \cup \Omega^{\ell}_0$, the
remainder of the ROMP algorithm is identical to OMP.

In order to analyze ROMP, we will need only two preliminary lemmas
from~\cite{romp}, which we state without proof.  Note that
Lemma~\ref{lem:prop32}, which is essentially a generalization of
Lemma~\ref{lem:ip2}, is stated using slightly weaker assumptions
than those stated in~\cite{romp}. The present version can easily be
obtained using the same proof.
%
\begin{lemma}[(1) in Prop. 3.2 of~\cite{romp}] \label{lem:prop32}
Let $\Gamma \subset \{1,2,\ldots,N\}$ and $x \in \reals^N$ be given.
Then if $\Psi$ satisfies the RIP of order  $| \supp(x) \cup \Gamma|$
with isometry constant $\delta$, we have
$$
\| (\Psi^T \Psi x) |_\Gamma - x|_\Gamma \|_2 \le \delta \|x\|_2.
$$
\end{lemma}
\begin{lemma} [Lemma 3.7 of~\cite{romp}] \label{lem:lem37}
Let $u \in \reals^K$, $K>1$, be arbitrary.  Then there exists a
subset $\Gamma \subseteq \{1, \ldots, K\}$ such that $|u(i)| \le 2
|u(j)|$ for all $i,j\in \Gamma$ and
$$
\|u|_\Gamma \|_2 \ge \frac{1}{2.5 \sqrt{\log_2 K } } \|u\|_2.
$$
\end{lemma}

Using these lemmas, we now provide a simplified proof of the main
result of~\cite{romp} concerning the recovery of sparse signals
using ROMP.\footnote{Note that we assume that $\Phi$ satisfies the
RIP of order $3K$ with constant $\delta \le 0.13/ \sqrt{\log_2 K}$.
Using Corollary 3.4 of~\cite{cosamp}, we can replace this with the
assumption that $\Phi$ satisfies the RIP of order $2K$ with constant
$\delta \le .043/ \sqrt{\log_2 K}$.}
\begin{thm}
Suppose that $\Phi$ satisfies the RIP of order $3K$ with isometry
constant $\delta \le 0.13/\sqrt{\log_2 K}$.
Then for any $x \in \reals^N$ with $\|x\|_0 \le K$, ROMP will
recover $x$ exactly from $y = \Phi x$ in at most $K$ iterations.
\end{thm}
~~~~~{\em Proof:} The proof works by showing that at each iteration,
\begin{equation} \label{eq:success}
| \Omega_0^\ell \cap \supp(x)| \ge \frac{1}{2} | \Omega_0^\ell |.
\end{equation}
If (\ref{eq:success}) is satisfied for $0,1,\dots,\ell-1$, then at
iteration $\ell$ we have that
\begin{equation} \label{eq:goodlambda}
| \Lambda^\ell \cap \supp(x)| \ge \frac{1}{2} |\Lambda^\ell |.
\end{equation}
It follows that, before $|\Lambda^\ell |$ exceeds $2K$, we will have
$\supp(x) \subseteq \Lambda^\ell$. Because $\Phi$ satisfies the RIP
of order $3K > 2K$, at termination, $\Phi_{\Lambda^\ell}$ will be
full rank. From (\ref{eq:ls}) we conclude that $x^{\ell} = x$
exactly.


To prove (\ref{eq:success}), we again proceed by induction. Hence,
we assume that (\ref{eq:success}) holds for $0, 1, \ldots, \ell-1$,
and thus (\ref{eq:goodlambda}) holds for iteration $\ell$.
%
%
We next assume for the sake of a contradiction that
(\ref{eq:success}) does not hold for iteration $\ell$, i.e., that
\begin{equation} \label{eq:cont}
| \Omega_0^\ell \setminus \supp(x)| > \frac{1}{2} |\Omega_0^\ell|.
\end{equation}
Define the sets $T = \Omega_0^\ell \setminus \supp(x)$ and $S =
\supp(x) \setminus \Lambda^\ell = \supp(\widetilde{x}^\ell)$, where
$\widetilde{x}^\ell$ is defined as in (\ref{eq:xtilde}).
Recall that we can write $h^\ell = A_{\Lambda^\ell}^T
A_{\Lambda^\ell} \widetilde{x}^\ell$.  Thus, using the assumption
that $|T| > \frac{1}{2} |\Omega_0^\ell|$ and the facts that $T
\subseteq \Omega_0^\ell$ and $\Omega_0^\ell \in R(\Omega^\ell)$,
one can show that
\begin{equation} \label{eq:romp1}
\| h^\ell |_T \|_2 \ge \frac{1}{\sqrt{5}} \|h^\ell |_{\Omega_0^\ell}
\|_2.
\end{equation}
We now observe that
\begin{equation} \label{eq:romp2}
\| h^\ell |_{\Omega_0^\ell} \|_2 \ge \frac{1}{2.5\sqrt{\log_2 K}} \| h^\ell|_{\Omega^\ell} \|_2,
\end{equation}
which follows from Lemma~\ref{lem:lem37} and the fact that $\Omega_0^\ell$ is the maximal regularizing set.  From the maximality of $\Omega^\ell$ and the fact that $|S| \le K$, we have that $\| h^\ell|_{\Omega^\ell} \|_2 \ge \| h^\ell |_S \|_2$, so that by combining (\ref{eq:romp1}) and (\ref{eq:romp2}) we obtain
\begin{equation} \label{eq:romp3}
\| h^\ell |_T \|_2 \ge \frac{1}{2.5 \sqrt{5 \log_2 K  }} \| h^\ell |_S \|_2.
\end{equation}
%
Note that $|S \cup \supp(\widetilde{x}^\ell)| = |S| \le K$ and since
$|\Lambda^\ell| \le 2K$, from Lemma~\ref{lem:PRIP}
we have that $A_{\Lambda^\ell}$ satisfies the RIP of order at least
$K$ with constant $\delta/(1-\delta)$, thus Lemma~\ref{lem:prop32}
implies that
$$
\|h^\ell|_S - \widetilde{x}^\ell|_S\|_2 \le \frac{\delta}{1-\delta} \|\widetilde{x}^\ell\|_2.
$$
Since $\widetilde{x}^\ell|_S = \widetilde{x}^\ell$,  $ \|h^\ell|_S -
\widetilde{x}^\ell|_S\|_2 = \|\widetilde{x}^\ell- h^\ell|_S\|_2 \ge
\|\widetilde{x}^\ell\|_2 - \|h^\ell|_S\|_2$,
and thus
$$
\|h^\ell|_S \|_2 \ge \frac{1-2\delta}{1-\delta} \|\widetilde{x}^\ell \|_2.
$$
Hence,
\begin{equation} \label{eq:romp4}
\| h^\ell |_T \|_2 \ge \frac{(1-2\delta)/(1-\delta)}{2.5 \sqrt{5 \log_2 K  }} \|\widetilde{x}^\ell\|_2.
\end{equation}
On the other hand, since $|\supp(\widetilde{x}^\ell)| +
|\Lambda^\ell \cap \supp(x)| = K$, from (\ref{eq:goodlambda}) we
obtain that $|\supp(\widetilde{x}^\ell)| \le K - \frac12
|\Lambda^\ell|$.  Thus,  $|T \cup \supp(\widetilde{x}^\ell)| \le |T|
+ |\supp(\widetilde{x}^\ell)| \le 2K - \frac12 |\Lambda^\ell|$.
Furthermore,  $A_{\Lambda^\ell}$ satisfies the RIP of order $3K -
|\Lambda^\ell| = 3K - \frac12 |\Lambda^\ell|- \frac12
|\Lambda^\ell|$.  Since $|\Lambda^\ell| \le 2K$, we have that
$A_{\Lambda^\ell}$ satisfies the RIP of order at least $2K - \frac12
|\Lambda^\ell|$ with constant $\delta/(1-\delta)$.
Thus, Lemma~\ref{lem:prop32} also implies that
\begin{equation} \label{eq:romp5}
\| h^\ell |_T \|_2 = \| h^\ell |_T - \widetilde{x}^\ell|_T \|_2 \le \frac{\delta}{1-\delta} \|\widetilde{x}^\ell\|_2.
\end{equation}
This is a contradiction whenever the right-hand-side of
(\ref{eq:romp4}) is greater than the right-hand-side of
(\ref{eq:romp5}), which occurs when $\delta < 1/(2 + 2.5 \sqrt{5
\log_2 K })$. Since $\log_2 K \ge 1$, we can replace this with the
slightly stricter condition $\delta < 1/((2+2.5\sqrt{5})\sqrt{\log_2
K}) \approx 0.1317/\sqrt{\log_2 K}$.
\qed \vspace{3mm}

Observe that when $K=1$, this proof (as well as the proofs
in~\cite{romp,romp2}) break down since Lemma~\ref{lem:lem37} does
not apply.  
However, when $K=1$ the ROMP algorithm simply reduces to OMP.  In this case we can apply Theorem~\ref{theo:main} to verify that ROMP succeeds when $K=1$ provided that $\Phi$ satisfies the RIP of order 2 with isometry constant $\delta < 1/3$.

\section*{Appendix}

\noindent{\em Proof of Lemma \ref{lem:ip}:}
We first assume that $\|u\|_2 = \|v\|_2 = 1$.  From the fact that
$$
\|u \pm v\|_2^2 = \|u\|_2^2 + \|v\|_2^2 \pm 2 \langle u,v\rangle = 2
\pm 2\langle u,v \rangle$$ and since $\Psi$ satisfies the RIP, we
have that
$$
(1-\delta)(2 \pm 2\langle u,v \rangle) \le \|\Psi u \pm \Psi v \|_2^2 \le (1+\delta)(2 \pm 2\langle u,v \rangle).
$$
From the parallelogram identity we obtain
$$
\langle \Psi u, \Psi v \rangle = \frac{1}{4} \left( \|\Psi u +
\Psi v \|_2^2 - \|\Psi u - \Psi v \|_2^2 \right) \le \frac{(1 +
\langle u,v \rangle)(1+\delta) - (1 -  \langle u,v
\rangle)(1-\delta)}{2}
 = \langle u,v \rangle + \delta.
$$
Similarly, one can show that $\langle \Psi u, \Psi v \rangle \ge
\langle u,v \rangle - \delta$, and thus
$|\langle \Psi u, \Psi v \rangle - \langle u,v \rangle| \le \delta$.  The result follows for
$u$, $v$ with arbitrary norm from the bilinearity of the inner product. \qed 

\vspace{3mm} \noindent{\em Proof of Lemma \ref{lem:PRIP}:}
From the definition of $A_\Lambda$ we may decompose $A_\Lambda u$ as
$A_\Lambda u = \Phi u - P_\Lambda \Phi u$. Since $P_\Lambda$ is an
orthogonal projection,  we can write
\begin{equation} \label{eq:triangle}
\|\Phi u\|_2^2 = \|P_\Lambda \Phi u\|_2^2 + \|A_\Lambda u\|_2^2.
\end{equation}
Our goal is to show that $\| \Phi u \|_2 \approx \|A_\Lambda u\|_2$,
or equivalently, that $\|P_\Lambda \Phi u\|_2$ is small. Towards
this end, we note that since $P_\Lambda \Phi u$ is orthogonal to
$A_\Lambda u$,
\begin{equation}
\langle P_\Lambda \Phi u, \Phi u \rangle  = \langle P_\Lambda \Phi
u, P_\Lambda \Phi u + A_\Lambda u \rangle
 = \langle P_\Lambda \Phi u, P_\Lambda \Phi u \rangle + \langle P_\Lambda \Phi u, A_\Lambda u \rangle
 = \|P_\Lambda \Phi u\|_2^2. \label{eq:cosine}
\end{equation}
Since $P_\Lambda$ is a projection onto $\Ran(\Phi_\Lambda)$ there
exists a $z \in \real^N$
with $\supp(z) \subseteq \Lambda$ such that $P_\Lambda \Phi u = \Phi
z$. Furthermore, by assumption, $\supp(u) \cap \Lambda = \emptyset$.
Hence $\langle u, z \rangle = 0$ and from the RIP and Lemma
\ref{lem:ip},
$$
\frac{|\langle P_\Lambda \Phi u, \Phi u \rangle|}{\|P_\Lambda \Phi
u\|_2\|\Phi u\|_2} = \frac{|\langle \Phi z, \Phi u \rangle|}{\| \Phi
z\|_2\|\Phi u\|_2}  \le  \frac{|\langle \Phi z, \Phi u \rangle|}{(1-\delta) \|z\|_2\|u\|_2} \le \frac{\delta}{1-\delta}.
$$
Combining this with (\ref{eq:cosine}), we obtain
$$
\| P_\Lambda \Phi u\|_2 \le \frac{\delta}{1-\delta} \|\Phi u\|_2.
$$
Since we trivially have that $\|P_\Lambda \Phi u\|_2 \ge 0$, we can
combine this with (\ref{eq:triangle}) to obtain
$$
\left(1- \left(\frac{\delta}{1-\delta}\right)^2\right)\|\Phi u\|_2^2
\le \| A_\Lambda u\|_2^2 \le \|\Phi u\|_2^2.
$$
Since $\|u\|_0 \le K$, we can use the RIP to obtain
$$
\left(1- \left(\frac{\delta}{1-\delta}\right)^2\right) (1-\delta)\|u\|_2^2 \le \| A_\Lambda u\|_2^2 \le (1+\delta)\|u\|_2^2,
$$
which simplifies to (\ref{eq:PRIP}). \qed 


\bibliographystyle{IEEEbib}
\footnotesize
\bibliography{OMPbib}


\end{document}